\documentclass[letterpaper, 10 pt, conference]{ieeeconf}
\IEEEoverridecommandlockouts   

\IEEEoverridecommandlockouts                              

\usepackage{graphicx}

\usepackage{url}
\usepackage{amsmath,amssymb,mathtools}
\usepackage{ thmtools}
\usepackage{amsfonts} 

\usepackage{enumerate}
\usepackage{color}                  

\usepackage{verbatim}

\usepackage{amssymb}
\usepackage{amsbsy}
\usepackage{amsmath}
\usepackage{setspace}
\usepackage{url}
\usepackage{subcaption}
\usepackage{float}
\usepackage{dsfont}

 \newcommand{\Int}{\operatorname{int}}

 \newcommand{\diag}{\operatorname{diag}}
 
 \newcommand{\tr}{\operatorname{tr}}

\newcommand{\trace}{\operatorname{trace}}
\newcommand\scalemath[2]{\scalebox{#1}{\mbox{\ensuremath{\displaystyle #2}}}}

\declaretheorem[name={Example}  ] {Example}

\declaretheorem[name={Definition}  ] {Definition}
\declaretheorem[name={Theorem}  ] {Theorem}
\declaretheorem[name={Lemma}  ] {Lemma}
\declaretheorem[name={Remark}  ] {Remark}

\declaretheorem[name={Proposition}  ] {Proposition}

\newcommand {\R}{\mathbb R}

\newcommand {\C}{\mathbb C}

\newcommand{\be}{\begin{equation}}
\newcommand{\ee}{\end{equation}}
\newcommand{\sgn}{\operatorname{{\mathrm sgn}}}

\newcommand*\dif{\mathop{}\!\mathrm{d}}

\usepackage{lineno}

\usepackage{todonotes}

 \title{\LARGE \bf
Compound matrices in systems and control theory }
 
\author{Eyal Bar-Shalom and   Michael Margaliot
\thanks{This research was  partially supported by a research grant from the Israel Science Foundation~(ISF)}
\thanks{The authors  are  with 
 the School of Electrical  Engineering,
		and the Sagol School of Neuroscience, 
		Tel-Aviv University, Tel-Aviv~69978, Israel.
		E-mail: \texttt{michaelm@tauex.tau.ac.il
		}}}
\begin{document}
\maketitle
%
\begin{abstract}
%
The  multiplicative and  additive compounds of a matrix play an important role in several fields of mathematics including geometry, 
multi-linear algebra, combinatorics, and the  analysis of
nonlinear time-varying dynamical  systems. There is a growing interest in applications of these compounds, and their generalizations, in systems and control theory. 
This tutorial paper provides a gentle introduction to these topics with an emphasis on  the geometric interpretation of the compounds, and surveys some of their recent applications. 
%
\end{abstract}
%
\section{Introduction}
%
 Let~$A\in\R^{n\times n}$. Fix~$k\in\{1,\dots,n\}$.
 The $k$-multiplicative and $k$-additive compounds of~$A$ play an important role in several fields of applied mathematics. Recently, there is a growing interest in the applications of these compounds, and their generalizations, in systems and control theory (see, e.g.~\cite{cheng_diag_stab,kordercont,rami_osci,rola_spect,CTPDS,margaliot2019revisiting,  pines2021,Eyal_k_posi,DT_K_POSI,9107214, gruss1,gruss2,gruss3,gruss4,gruss5}). 
 
 This tutorial paper  reviews the~$k$-compounds,  focusing on their  geometric interpretations, and
 surveys some of their  recent applications in systems and control theory,  including~$k$-positive systems,
 $k$-cooperative systems, 
 and~$k$-contracting systems. 
 
   The results described here are known, albeit never collected in a single paper. The only exception is the 
new    generalization 
principle  described in Section~\ref{sec:k_generalizations}.

We use standard notation. For a set~$S$, $\Int(S)$ is the interior of~$S$.
For scalars~$\lambda_i$, $i\in\{1,\dots, n\}$, $\diag(\lambda_1,\dots,\lambda_n)$ is the $n\times n$ diagonal matrix with
diagonal entries~$\lambda_i$. 
Column vectors [matrices] are denoted by small [capital] letters. For a matrix~$A$, $A^T$ is the transpose of~$A$. For a square matrix~$B$, $\tr(B)$ [$\det(B)$] is the trace [determinant]  of~$B$. 
 $B $ is called Metzler if all its off-diagonal entries are non-negative.  The positive orthant in~$\R^n$ is~$\R^n_+ = \{x\in\R^n: x_i\geq 0,\; i=1,\dots,n\}$.  

\section{Geometric motivation}
\begin{figure}
 \begin{center}
  \includegraphics[scale=0.33]{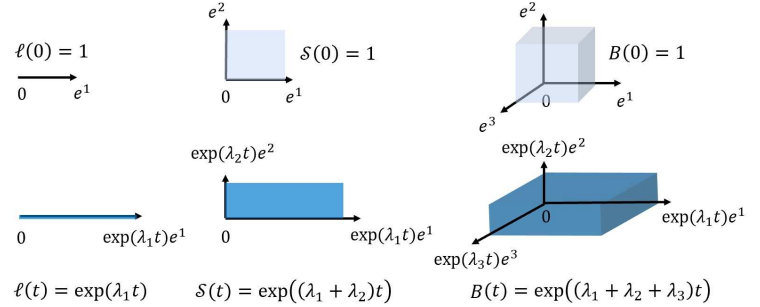}
	\caption{The evolution of lines, areas, and volumes  under the   the LTI~\eqref{eq:ltia}.}
	\label{fig:GeometricalEvolution}
\end{center}
\end{figure}

$k$-compound matrices provide information on the evolution of $k$-dimensional polygons   subject  to    a linear time-varying dynamics. To explain this in a simple setting,
consider the LTI system:
\be\label{eq:ltia}
\dot x =\diag(\lambda_1,\lambda_2,\lambda_3)x ,
\ee
with~$\lambda_i\in\R$ and~$x:\R_+ \to\R^3$.
Let~$e^i$, $i=1,2,3$,
denote the~$i$th canonical vector in~$\R^3$.
For~$x(0)=e^i$ we have~$x(t)=\exp(\lambda_i t)x(0)$. Thus, $\exp(\lambda_i t)$ describes the rate of evolution of the line~$e^i$ subject to~\eqref{eq:ltia}.
What about 2D areas? Let~$S\subset\R^3 $ denote the square generated by~$e^i$ and~$e^j$, with~$i\not =j$.
Then~$S(t):= x(t,S )$ is the rectangle generated by~$\exp(\lambda_i t) e^i$ and~$\exp(\lambda_j t) e^j$, so 
  the area of~$S(t)$ is~$\exp( (\lambda_i + \lambda_j) t)$. Similarly, if~$B \subset\R^3 $ is the 3D box generated by~$e^1,e^2$, and~$e^3$ then the volume of~$B(t):=x(t,B)$  
is~$\exp( (\lambda_1+\lambda_2+\lambda_3)t)$ (see Fig.~\ref{fig:GeometricalEvolution}). 

Since~$\exp(At)=\diag ( \exp(\lambda_1 t), \exp(\lambda_2 t), \exp(\lambda_3 t)  )$, this discussion suggests
that it may be useful to have a~$3\times 3$   matrix whose eigenvalues are the sums of any two eigenvalues of~$\exp(At)$, 
and a~$1\times 1$ matrix whose eigenvalue
is the sum of the  three  eigenvalues of~$\exp(At)$.
With this geometric motivation in mind, we   turn to recall the notions
of the multiplicative and additive compounds of a matrix. For more details and proofs, see e.g.~\cite[Ch.~6]{fiedler_book}\cite{schwarz1970}.

%

\section{Compound Matrices}
%
Let~$A\in\C^{n\times m }$.   Fix~$k\in\{1,\dots,\min\{n,m\} \}$.
Let~$Q(k,n)$ denote the set of 
increasing sequences of~$k$ integers in~$\{1,\dots,n\}$, ordered lexicographically.
For example,
\[
Q(2,3)=\{ \{1,2\} ,\{1,3\},\{2,3\} \}.
\]
For~$\alpha  \in Q(k,n),\beta \in Q(k,m)$, let~$A[\alpha|\beta]$ denote the $ k\times k$ submtarix   obtained by taking the entries of~$A$ in  the rows indexed by~$\alpha$ and the columns indexed by~$\beta$. For example
\[
A[ \{1,2\}|\{2,3\}]=\begin{bmatrix}
a_{12} & a_{13}\\
a_{22} &a_{23}
\end{bmatrix} .
\]
The minor of~$A $ corresponding to~$\alpha ,  \beta$ is~$A(\alpha|\beta):=\det (A[\alpha|\beta]) $. For example,~$Q(n,n)$ includes the single set~$\alpha=\{1,\dots,n\}$
and~$A(\alpha|\alpha)=\det(A)$.
%
\begin{Definition}\label{def:multi} 
%
Let $A\in C^{n\times m}$ 
and fix $k \in \{ 1,\dots,\min \{ n,m \} \}$. 
The \emph{$k$-multiplicative compound} of~$A$, denoted~$A^{(k)}$, is the~$\binom{n}{k}\times \binom{m}{k}$
matrix that contains all the~$k$-order minors of~$A$  ordered lexicographically. 
\end{Definition}

For example, if~$n=m=3$ and~$k=2$ then
\[
A^{(2)}= \begin{bmatrix}
A(\{12\}|\{12\}) & A(\{12\}|\{13\}) & A(\{12\}|\{23\}) \\ 
A(\{13\}|\{12\}) & A(\{13\}|\{13\}) & A(\{13\}|\{23\}) \\ 
A(\{23\}|\{12\}) & A(\{23\}|\{13\}) & A(\{23\}|\{23\}) 
\end{bmatrix}.
\]
In particular, Definition~\ref{def:multi} 
implies that~$A^{(1)}=A$,
and if~$n=m$ then~$A^{(n)}=\det(A)$.
%
 
%
%
%
Let~$A\in\mathbb{C}^{n\times m},\; B\in\mathbb{C}^{m\times p}$. The Cauchy-Binet Formula (see e.g.~\cite{notes_comb_algebra})
asserts 
that for
any   $k\in \{1,\dots,\min \{n,m,p\} \}$,
\begin{align}\label{eq:AB_MultComp}
    (AB)^{(k)} = A^{(k)} B^{(k)}.
\end{align}
Hence the term multiplicative compound.
%
%
Note that for~$n=m=p$,  Eq.~\eqref{eq:AB_MultComp}  with~$k=n$ reduces to the familiar formula~$\det(AB)= \det (A) \det (B)$.
%

Let~$I_s$ denote the~$s\times s$ identity matrix. 
Definition~\ref{def:multi}  implies that~$I_n^{(k)} = I_r$, where~$r:=\binom{n}{k}$. Hence, if~$A\in\R^{n\times n}$ is non-singular then~$(AA^{-1})^{(k)}=I_r$ and combining this with~\eqref{eq:AB_MultComp} yields~$(A^{-1})^{(k)} = (A^{(k)})^{-1}$.
In particular, if~$A$ is non-singular then so is~$A^{(k)}$. 
 
%
 
%

The~$k$-multiplicative compound has an important spectral property. 
For~$A\in\C^{n\times n}$, let~$\lambda_i$,~$i=1,\dots,n$, denote the eigenvalues of~$A$. 
Then the eigenvalues of~$A^{(k)}$ are all the products 
\be\label{eq:prodeig}
\prod_{\ell=1}^k \lambda_{i_\ell},
\text{
with } 1\leq i_1< i_2<\dots< i_k \leq n.
\ee

For example, suppose that~$n=3$ and  
$
A=\begin{bmatrix}
a_{11} & a_{12} & a_{13}\\
0 & a_{22} & a_{23}\\
 0 &  0  & a_{ 33}
\end{bmatrix}  .
$
Then a  calculation gives
\[
A^{(2)}= \begin{bmatrix}
a_{11}   a_{22} &  a_{11} a_{23}&  a_{12} a_{23}-a_{13} a_{22} \\
0& a_{11} a_{33}& a_{12} a_{33}\\
0& 0& a_{22} a_{33}
 \end{bmatrix} ,
\]
so, clearly, the eigenvalues of~$A^{(2)}$
are of the  form~\eqref{eq:prodeig}. 
\begin{Definition}
%
Let~$A\in\mathbb{C}^{n\times n}$.
The \emph{$k$-additive compound}
matrix of~$A$ is defined by:
\begin{align*}
    A^{[k]} := \frac{d}{d\epsilon} (I+\epsilon A)^{(k)} |_{\epsilon=0} . 
\end{align*}
%
\end{Definition}

Note that this implies that~$ A^{[k]} =  \frac{d}{d\epsilon} (\exp(A\epsilon ))^{(k)} |_{\epsilon=0}$, and also that
\begin{align}\label{eq:(I+epsilonA)^k}
    (I+\epsilon A)^{(k)} = I + \epsilon A^{[k]} + o(\epsilon),
\end{align}
%
 In other words,~$A^{[k]}$ is the first-order term in the Taylor expansion of~$(I+\epsilon A)^{(k)}$.
 
Let~$\lambda_i$,~$i=1,\dots,n$, denote the eigenvalues of~$A$. Then the eigenvalues of~$(I+\epsilon A)^{(k)}$ are the products~$\prod_{\ell=1}^k (1+\epsilon \lambda_{i_\ell}) $, and~\eqref{eq:(I+epsilonA)^k} implies that 
  the eigenvalues of~$A^{[k]}$ are all the sums 
\[
\sum_{\ell=1}^k \lambda_{i_\ell},
\text{
with } 1\leq i_1< i_2<\dots< i_k \leq n.
\]

Another important  implication  of the definitions above is that for any~$A,B\in\C^{n\times n}$ we have
\[
(A+B)^{[k]}=A^{[k]}+B^{[k]}.
\]
This justifies the term additive compound. 
Moreover, the mapping~$A\to A^{[k]}$ is linear.

The following result gives a useful explicit formula for~$A^{[k]}$ in terms of the entries~$a_{ij}$ of~$A$.
  Recall that any entry of~$A^{(k)}$ is a minor~$A(\alpha|\beta)$.
Thus, it is natural to index the entries of~$A^{(k)}$  and~$A^{[k]}$  using~$\alpha,\beta \in Q(k,n)$. 
%
\begin{Proposition}\label{prop:Explicit_A_k}
%
Fix~$\alpha,\beta \in Q(k,n)$ and let~$\alpha=\{i_1,\dots,i_k\}$ and~$\beta=\{j_1,\dots,j_k\}$. Then the entry of~$A^{[k]}$ corresponding to~$(\alpha,\beta)$ is
equal to:
\begin{enumerate}
    \item $\sum_{\ell=1}^{k} a_{ i_{\ell} i_{\ell} }$, if  $i_{\ell} = j_{\ell}$  for all $\ell \in \{ 1,\hdots,k \}$;
    \item $(-1)^{\ell +m}
    a_{i_{\ell} j_{m}} $, 
      if all the indices  in $ \alpha  $  and $ \beta$   agree,  except  for   
   a single index $ i_{\ell} \ne j_m$; and
    \item $0$, otherwise.
\end{enumerate}
%
\end{Proposition}

\begin{Example}\label{ex:A^[3]}
%
For~$A\in\R^{4\times 4}$ and~$k=3$,
Prop.~\ref{prop:Explicit_A_k} yields
\[ 
A^{[3]}  =
\left[ \scalemath{0.67}{
\begin{array}{ccccccc} 
    a_{11}+a_{22}+a_{33} & a_{34} & -a_{24} & a_{14} \\
    a_{43} & a_{11}+a_{22}+a_{44} & a_{23} & -a_{13} \\
    -a_{42} & a_{32} & a_{11}+a_{33}+a_{44} & a_{12} \\
    a_{41} & -a_{31} & a_{21} & a_{22}+a_{33}+a_{44}
\end{array}} \right].
\]
The entry in the second row and fourth column
of $A^{[3]}$ corresponds to 
$(\alpha | \beta) = ( \{ 1,2,4 \}  |  \{ 2,3,4 \} )$.
As $\alpha$ and $\beta$ agree in all indices except
for the  index
$i_{ 1} =1$ and $j_{ 2} =3$,
this  entry is equal to~$  (-1)^{1+2} a_{13} = -a_{13}$.
%
\end{Example}

Note that
Prop.~\ref{prop:Explicit_A_k}
implies in particular that~$A^{[n]}=\tr(A)$.

The next section describes applications of compound matrices for dynamical systems described by~ODEs. For more details and proofs, see~\cite{schwarz1970,muldowney1990}.
\section{Compound Matrices and ODEs}
%
Fix an interval~$[a,b] \in \R$. Let~$A:[a,b] \to \R^{n\times n} $ be a continuous matrix function, and consider the LTV system:
\be\label{Eq:ltv}
\dot x(t)=A(t)x(t),\quad x(a)=x_0. 
\ee
The solution is given by~$x(t)=\Phi(t ,a)x(a)$,
where~$\Phi(t,a )$ is the solution at time~$t$ of the matrix differential 
equation
\be\label{eq:phidot}
\dot \Phi(s)=A(s) \Phi(s), 
\quad \Phi(a)=I_n. 
\ee
Fix~$k\in\{1,\dots,n\}$ and let~$r:=\binom{n}{k}$.
A natural question is: how do the~$k$-order minors of~$\Phi(t)$ evolve in time?
The next result provides a beautiful formula for the evolution of~$\Phi^{(k)}(t):= (\Phi(t))^{(k)}$.
\begin{Proposition}\label{prop:odeforphik}
If~$\Phi$ satisfies~\eqref{eq:phidot} then
\be\label{eq:expat}
\frac{d}{dt} \Phi^{(k)}(t)=A^{[k]}(t) \Phi^{(k)}(t),\quad \Phi^{(k)}(a)=I_r,
\ee
where~$A^{[k]}(t):= (A(t)) ^ {[k]} $.
\end{Proposition}

Thus, the~$k\times k$ minors of~$\Phi$ also satisfy an~LTV.
In particular, if~$A(t)\equiv A$  and~$a=0$
then~$\Phi(t)=\exp( At) $ so~$\Phi^{(k)}(t)=(\exp(At))^{(k)}$ and~\eqref{eq:expat} gives
\[
(\exp(At))^{(k)} = 
\exp(A ^{[k]} t).
\]

Note also that for~$k=n$, Prop.~\ref{prop:odeforphik}  yields 
 \[
 \frac{d}{dt} \det(\Phi(t))
 = \tr(A(t)) \det(\Phi(t)),
 \]
which is the
Abel-Jacobi-Liouville identity.

Roughly speaking, Prop.~\ref{prop:odeforphik} implies that under the LTV dynamics~\eqref{Eq:ltv}, $k$-dimensional polygons evolve  according to the dynamics~\eqref{eq:expat}.

We now turn to consider the nonlinear system:
\be\label{eq:nonlin}
\dot x (t) = f(t,x).
\ee
For the sake of simplicity, we  assume that the initial time is zero, and that the system admits a convex and compact state-space~$\Omega$. We also assume that~$f\in C^1$. The Jacobian of the vector field~$f$ is~$J(t,x):=\frac{\partial}{\partial x}f(t,x)$.

Compound matrices can be used to analyze~\eqref{eq:nonlin} by using an LTV called the  variational equation associated  with~\eqref{eq:nonlin}.
To define it, 
fix~$a,b \in \Omega$.
Let~$z(t):=x(t,a)-x(t,b)$, and for~$s\in[0,1]$, 
 let~$\gamma(s):=s x(t,a)+(1-s)x(t,b)$, i.e. the line connecting~$x(t,a)$ and~$x(t,b)$.
 Then
\begin{align*} 
   \dot z(t)&=f(t,x(t,a))-f(t,x(t,b))  \\
	       &=\int_0^1  \frac{\partial }{\partial s}   f(t,\gamma(s))     \dif s  ,
\end{align*}
and this gives the variational equation:
\be\label{eq:var_eqn}
\dot z(t)=A^{ab}(t)z(t),
\ee
where
\begin{equation}\label{eq:at_int}
    A^{ab}(t):=\int_0^1  J(t,\gamma(s))     \dif s.  
\end{equation}

We can use the results above to  describe a powerful approach  for deriving useful
``$k$-generalizations'' of  important 
classes of dynamical systems including cooperative systems~\cite{hlsmith},  contracting systems~\cite{sontag_cotraction_tutorial,LOHMILLER1998683}, 
and diagonally stable systems~\cite{diag_Stab_book}. 

\section{$k$-generalizations of dynamical systems}\label{sec:k_generalizations}
Consider the LTV~\eqref{Eq:ltv}. Suppose that~$A(t)$ satisfies a specific \emph{property}, e.g. \emph{property} may be that~$A(t)$ is Metzler for all~$t$ (so the LTV is positive) or that~$\mu(A(t))\leq-\eta<0$ for all~$t$, where~$\mu:\R^{n\times n}\to\R$ is a matrix measure (so the~LTV is contracting). Fix~$k\in\{1,\dots,n\}$.  We say that
the LTV satisfies~\emph{$k$-property} if~$A^{[k]}$ (rather than~$A$)
satisfies \emph{property}. For example, the LTV is~\emph{$k$-positive} if
$A^{[k]}(t)$ is Metzler for all~$t$;  the LTV is~\emph{$k$-contracting} if~$\mu(A^{[k]}(t))\leq-\eta<0$ for all~$t$, and so on. 

This generalization approach  makes sense for two reasons. First, when~$k=1$,~$A^{[k]}$ reduces to~$A^{[1]}=A$, so 
\emph{$k$-property} is clearly a generalization of
\emph{property}. Second, we know that~$A^{[k]}$ has a clear geometric meaning, as it describes the evolution of~$k$-dimensional polygons along the dynamics. 

The same idea can be applied to  the nonlinear system~\eqref{eq:nonlin}
using the variational equation~\eqref{eq:var_eqn}.
For example, if~$\mu(J(t,x))\leq -\eta<0$ for all~$t\geq 0$ and~$x\in\Omega$ then~\eqref{eq:nonlin} is contracting: the distance between any two solutions (in the norm that induced~$\mu$) decays at an exponential rate. If we replace this by the condition~$\mu(J^{[k]}(t,x))\leq -\eta<0$ for all~$t\geq 0$ and~$x\in\Omega$ then~\eqref{eq:nonlin} is called~$k$-contracting. Roughly speaking, 
this means that the  volume of $k$-dimensional 
polygons  decays to zero exponentially along the flow of the nonlinear system.
We now turn to describe several  such~$k$-generalizations.

\section{$k$-contracting   systems}
$k$-contracting systems were introduced in~\cite{kordercont}
(see also the unpublished preprint~\cite{weak_manchester} 
for some preliminary ideas).
For~$k=1$ these reduce to contracting systems. 
This generalization was motivated in part by the seminal work of Muldowney~\cite{muldowney1990}
who considered nonlinear systems that, in the new terminology, are~$2$-contracting. He derived several interesting results for time-invariant $2$-contracting systems. For example, every bounded trajectory of a time-invariant, nonlinear,  $2$-contracting system
converges to an equilibrium (but, unlike in the case of contracting systems, the equilibrium is not necessarily unique). 

For the sake of simplicity, we introduce the ideas in the context of an LTV system.
The analysis of nonlinear systems is based on assuming that their variational equation
is  a~$k$-contracting LTV. 

Recall that a vector norm $|\cdot |:\mathbb{R}^{n}\to\mathbb{R}_{+}$ induces a
matrix norm $||\cdot ||:\mathbb{R}^{n\times n}\to\mathbb{R}_{+}$   by: 
\begin{align*}
    ||A||:=\max_{|x|=1} |Ax|,
\end{align*}

and a matrix measure   
$\mu (\cdot ):\mathbb{R}^{n\times n}\to\mathbb{R}$  
  by:
\begin{align*}
    \mu (A) :=\lim_{\epsilon \downarrow 0} \frac{||I+\epsilon A||-1}{\epsilon}.
\end{align*}

For~$p\in\{1,2,\infty\}$, let~$|\cdot|_p:\R^n\to\R_+$ denote the~$L_p$ vector norm, that is,
 $|x|_1:=\sum_{i=1}^n|x_i|$, 
$|x|_2:=\sqrt{x^T x}$, 
and~$|x|_\infty:=\max_{   i\in\{1,\dots,n\}    } |x_i|$.  
The induced  matrix measures are~\cite{vid}:
\begin{equation}\label{eq:matirx_measures_12infty} \begin{aligned} 
\mu_1(A) &= \max_{j\in\{1,\dots,n\}}  (a_{jj}+ \sum_{\substack{i=1 \\ i \neq j}}^n |a_{ij}|   ) , \\
\mu_2(A) &= \lambda_1 (  {A + A^T} )/2 , \\
\mu_{\infty}(A) &= \max_{i\in\{1,\dots,n\}}  (a_{ii}+ \sum_{\substack{j=1 \\ j \neq i}}^n |a_{ij}|  ) , 
\end{aligned} \end{equation}
where~$\lambda_i (S)$ denotes the $i$-th largest eigenvalue of the symmetric matrix~$S$, that is,
\[
\lambda_1(S) \geq \lambda_2(S) \geq \cdots \geq \lambda_n(S).
\]

\begin{Definition}
The LTV~\eqref{Eq:ltv} is called~\emph{$k$-contracting} with respect to (w.r.t.) the norm~$|\cdot| $ 
if  
\be\label{eq:kconinfi} 
\mu(A^{[k]}(t))\leq-\eta<0 \text{ for all } t \geq 0,
\ee
where~$\mu$ is the matrix measure induced by~$|\cdot|$.
\end{Definition}

Note that for~$k=1$ condition~\eqref{eq:kconinfi}  reduces to the standard infinitesimal condition for contraction~\cite{sontag_cotraction_tutorial}. For the~$L_p$ norms, with~$p\in\{1,2,\infty\}$, 
this condition is easy to check using the 
  explicit expressions for~$\mu_1$, $\mu_2$, and~$\mu_\infty$. This carries over to~$k$-contraction, as 
  combining Prop.~\ref{prop:Explicit_A_k} with~\eqref{eq:matirx_measures_12infty} gives~\cite{muldowney1990}:
\begin{align*} 
\mu_1(A^{[k]}) &= 
 \max_{\alpha \in Q(k,n) }  ( \sum_{p=1}^k a_{\alpha_p,\alpha_p}  
 + \sum_{\substack{j \notin \alpha }}(|a_{j,\alpha_1}| + \cdots + |a_{j,\alpha_k}|) ) ,\nonumber \\
\mu_2(A^{[k]}) &= \sum_{i=1}^k \lambda_i (  {A + A^T}   )/2 , \\
\mu_{\infty}(A^{[k]}) &=
\max_{\alpha \in Q(k,n) } (  \sum_{p=1}^k a_{\alpha_p,\alpha_p} + \sum_{\substack{j \notin\alpha }}(|a_{\alpha_1,j}| + \cdots + |a_{\alpha_k,j}|) ) .
\end{align*}

For~$k=n$, $A^{[n]}$ is the scalar~$\tr(A)$, so 
condition~\eqref{eq:kconinfi} becomes~$\tr(A(t))\leq -\eta<0 $ for all $t \geq 0$.

Combining 
Coppel's inequality~\cite{coppel1965stability}  with~\eqref{eq:expat} yields  the following result.
\begin{Proposition}
If the LTV~\eqref{Eq:ltv} is~$k$-contracting then 
\[
||\Phi^{(k)}(t) || \leq \exp(-\eta t) || \Phi^{(k)}(0)|| = \exp(-\eta t) 
\]
for all~$t\geq 0$.
\end{Proposition}

Geometrically, this means that under the LTV dynamics 
the volume of $k$-dimensional polygons converges to zero exponentially. The next example illustrates this.  

\begin{figure}
 \begin{center}
  \includegraphics[scale=0.5]{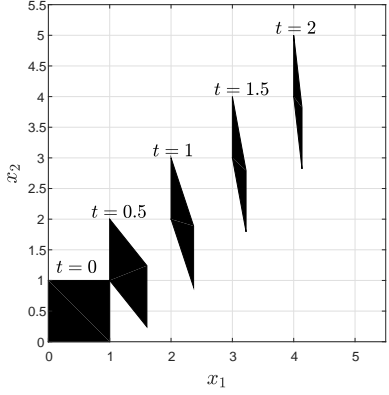}
	\caption{Evolution of  the unit square in Example~\ref{exa:squares}.}
	\label{fig:squares}
\end{center}
\end{figure}

\begin{Example}\label{exa:squares}
	Consider the LTV~\eqref{Eq:ltv} with~$n=2$ and
$
	A(t) = \begin{bmatrix}
	 -1 & 0 \\
	-2\cos( t) &0  
	\end{bmatrix}.
$
The corresponding  transition matrix is:
$
						  	\Phi(t)= \begin{bmatrix} \exp(-t)&0 \\ 
										                -1+\exp(-t)(\cos(t)-\sin(t))  &1\end{bmatrix}.
	$	This implies that the LTV is uniformly stable, and that for any~$x(0)\in\R^2$ we have
	\[
	\lim _{t\to \infty} x(t,x(0))=\begin{bmatrix}  0  \\ x_2(0)-x_1(0)  \end{bmatrix} .
	\]
	The LTV 
 is not contracting,  w.r.t.   any norm, as   there is more than a single equilibrium. However,~$ A^{[2]}(t) =\tr(A(t)) \equiv -1$, so the 
	system is~$2$-contracting.  Let~$S\subset\R^2$ denote the unit square, and let~$S(t):=x(t,S )$, that is,
	the evolution
	at time~$t$
	of the unit square   under the dynamics. 
	Fig.~\ref{fig:squares} 
	depicts~$S(t) +2t$  for several values of~$t$,
	where the shift by~$2t$ is for  the sake of clarity.
	It may be seen that the area of~$S(t)$ decays with~$t$, and that~$S(t)$ converges to a line.
\end{Example}

As noted above,   time-invariant $2$-contracting systems have a ``well-odered''
asymptotic behaviour~\cite{muldowney1990,li1995}, and this has been used to derive a global analysis of important models from epidemiology (see, e.g.~\cite{SEIR_LI_MULD1995}). 
A recent paper~\cite{searial12contracting} extended some of these results to systems that are not necessarily~$2$-contracting, but
can be represented as the serial interconnections of~$k$-contracting systems, with~$k\in\{1,2\}$. 

\section{$\alpha$-compounds and $\alpha$-contracting systems}
A  recent paper~\cite{pines2021} defined a generalizations
called the~$\alpha$-multiplicative compound and~$\alpha$-additive compound of a matrix, where~$\alpha$ is a \emph{real} number. Let~$A\in\C^{n\times n} $ be non-singular. 
If~$\alpha=k+s$, where~$k\in\{1,\dots,n-1\}$ and~$s\in(0,1)$ then the~$\alpha$-multiplicative of~$A$ is defined by:
\[
A^{(\alpha)} : = (A^{(k)})^{1-s} \otimes (A^{(k+1)})^{s},
\]
where~$\otimes$ denotes the Kronecker product. This is a kind of ``multiplicative interpolation'' between~$ A^{(k)} $ and~$A^{(k+1)} $.
For example,~$A^{(2.1)}    = (A^{(2)})^{0.9} \otimes (A^{(3)})^{0.1}$. The~$\alpha$-additive compound is defined just like the~$k$-additive compound, that is, 
\begin{align*}
    A^{[\alpha]} := \frac{d}{d\epsilon} (I+\epsilon A)^{(\alpha)} |_{\epsilon=0} , 
\end{align*}
and it was shown in~\cite{pines2021}
that this yields
\[
 A^{[\alpha]} = ((1-s)A^{[k]}) \oplus (sA^{[k+1]}),
 \]
where~$\oplus $ denotes  the Kronecker sum. 

The system~\eqref{eq:nonlin}
is called \emph{$\alpha$-contracting} w.r.t. the norm~$|\cdot| $ if
\be \label{eq:alphacon}
\mu(J^{[\alpha]}(t,x))\leq -\eta<0 ,
\ee
for all~$t\geq 0$ and~all~$x$ in the state space~\cite{pines2021}. 

Using this notion, it is possible to restate the seminal results of  Douady 
 and Oesterl\'{e}~\cite{Douady1980}  as follows. 
\begin{Theorem}\cite{pines2021}
Suppose that~\eqref{eq:nonlin}
is~$\alpha$-contracting for some~$\alpha\in [1,n]$.
Then any compact and strongly invariant set of the dynamics 
has a  Hausdorff dimension smaller than~$\alpha$. 
\end{Theorem}

Roughly speaking, the dynamics contracts   sets with
a  larger Hausdorff dimension.

\begin{figure}[t]
 \begin{center}
\includegraphics[width=8cm,height=6cm]{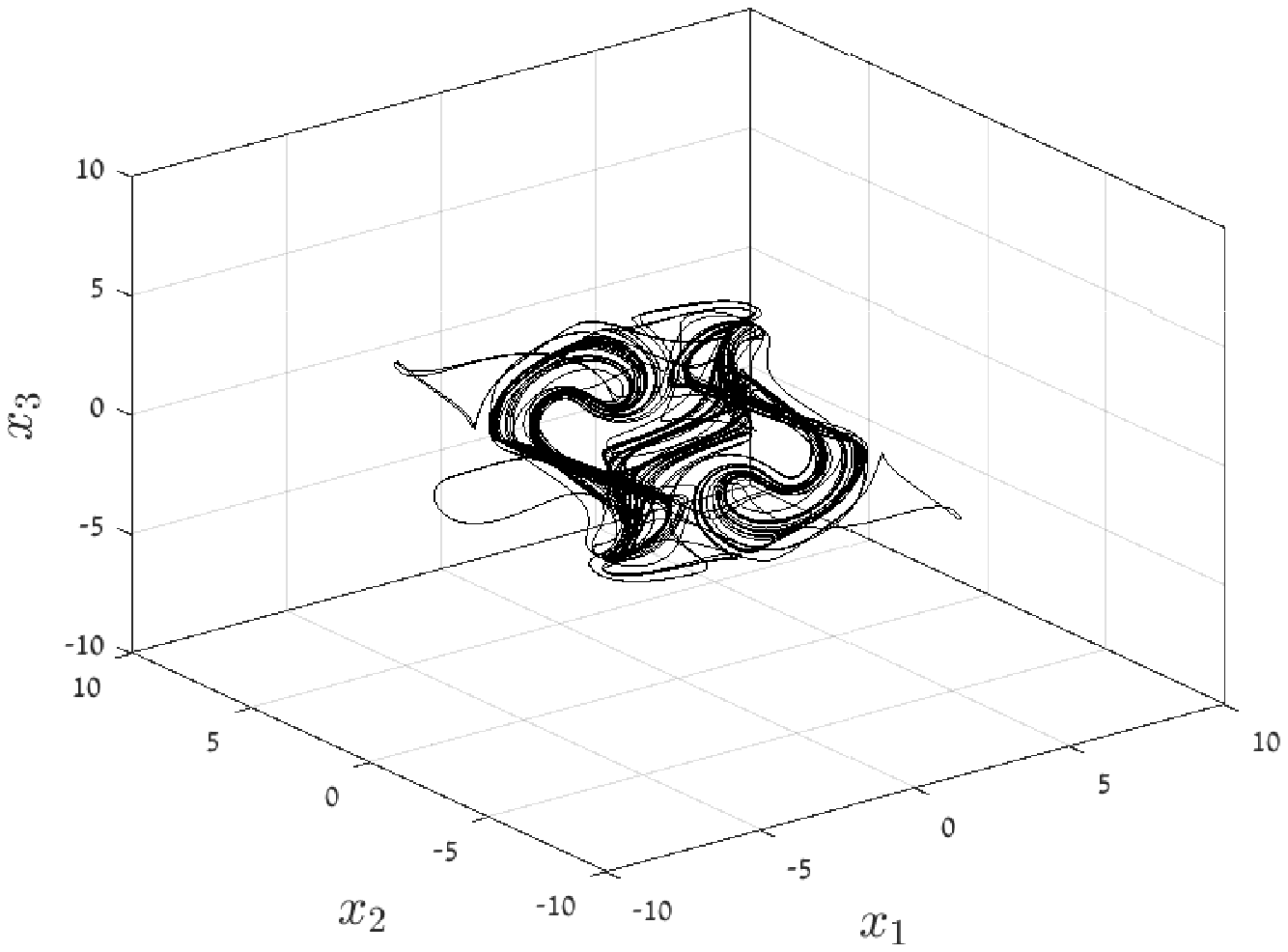}
\caption{A trajectory  emanating from~$x(0)=\begin{bmatrix} 
  1& -2&1 \end{bmatrix}^T$.
}\label{fig:chaos}
\end{center}
\end{figure}

\begin{figure}[t]
 \begin{center}
\includegraphics[width=8cm,height=6cm]{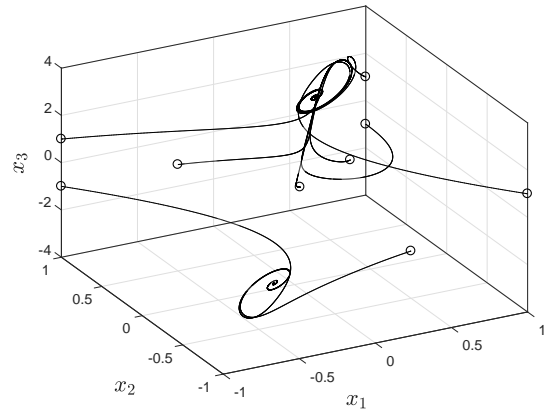}
\caption{Several trajectories of  the closed-loop system.  The circles denote the initial conditions of the trajectories.  }\label{fig:chaos_closed}
\end{center}
\end{figure}

The next example, adapted from~\cite{pines2021}, shows how these notions  can be used to ``de-chaotify'' a nonlinear dynamical system by feedback. 
\begin{Example}
Thomas' cyclically symmetric attractor~\cite{thomas99,chaos_survey} is
a popular example for a chaotic system. It is described by:
\begin{align} \label{eq:thom}
\dot x_1 =&  \sin(x_2)-bx_1 , \nonumber \\
\dot x_2 =&  \sin(x_3) - bx_2, \\
\dot x_3 =& \sin(x_1) - bx_3, \nonumber
\end{align}
where $b>0$ is the
dissipation  constant. 
The convex and compact set~$D: = \{x\in\R^3: b |x|_\infty \leq 1 \}$ is an
invariant set of the dynamics.
 
 Fig.~\ref{fig:chaos}
 depicts  the solution of the system emanating from~$\begin{bmatrix} 1& -2&1\end{bmatrix}^T$ for
 $
 b=0.1$.
 Note the symmetric strange attractor. 
 
The Jacobian $J_f$ of the vector field in~\eqref{eq:thom} is
\begin{align*}
   J_f (x)=\begin{bmatrix}-b&\cos(x_2)&0 \\ 0&-b&\cos(x_3) \\ \cos(x_1)&0&-b\end{bmatrix},
   \end{align*}
   and thus
 \begin{align*}
   J_f ^{[2]}(x)=\begin{bmatrix}-2b&\cos(x_3)&0 \\ 0&-2b&\cos(x_2) \\ -\cos(x_1)&0&-2b\end{bmatrix},
   \end{align*}
and~$
   J_f ^{[3]}=\trace(J (x))=-3b$. 
Since~$b>0$, this implies that the system is $3$-contracting  w.r.t. any norm. Let~$\alpha = 2+s$, with~$s\in(0,1)$. Then
\begin{align*}
    &J_f ^{[\alpha]}(x)
    =(1-s)J_f^{[2]}(x)\oplus sJ_f^{[3]}(x)\\
    &=
    \begin{bmatrix}
    -(2+s)b & (1-s)\cos(x_3) & 0 \\
    0 & -(2+s)b & (1-s)\cos(x_2) \\
    -(1-s)\cos(x_1) & 0 & -(2+s)b \\
    \end{bmatrix}.
\end{align*}
This implies that
\be\nonumber
\mu_{1 }(J_f ^{[\alpha]}(x)) \leq 1-2b-s(b+1), \text{ for all } x\in D.
\ee
We conclude that for any~$b\in(0,1/2) $
 the system is~$(2+s)$-contracting for any
$
    s>\frac{1-2b}{1+b}
$.

We now show how $\alpha$-contarction can be used to design a partial-state   controller for the system guaranteeing that the closed-loop system has a ``well-ordered'' behaviour. 
Suppose that 
the closed-loop system is:
\[
\dot x = f(x) + g(x),
\]
where~$g $ is the controller.
Let~$\alpha = 2 + s$, with~$s \in(0,1)$. 
The Jacobian of the closed-loop system is~$J_{cl}:=J_f+J_g$,  so
\begin{align*}
    \mu_1(J_{cl}^{[\alpha]})& = \mu_1(J_f^{[\alpha]}+J_g^{[\alpha]}) \\& \leq \mu_1(J_f^{[\alpha]})+\mu_1(J_g^{[\alpha]}) \\
    &\leq 1-2b-s(b+1)+\mu_1(J_g^{[\alpha]}).
\end{align*}
This implies that the closed-loop system is~$\alpha$-contracting if
\begin{align}\label{eq:cond_cont}
    \mu_1(J^{[\alpha]}_g (x) ) <    s(b+1)+2b-1
    \text{ for all } x\in D . 
\end{align}
Consider, for example,   the controller  
$
g(x_1,x_2)=c \diag(1,1,0) x $, with gain~$c<0$.
Then
$
J_g^{[\alpha]} = c\diag(  2  ,1+ s,1+ s )
$
and for any~$c<0$  condition~\eqref{eq:cond_cont} becomes
\begin{align}\label{eq:cond_cont1}
     (1+ s)c <    s(b+1)+2b-1.
\end{align}
This provides a simple recipe  for 
determining the gain~$c$ so that the closed-loop system is~$(2+s)$-contracting. For example,  
  when~$s \to 0$, Eq.~\eqref{eq:cond_cont1} yields
$
c<2b-1
$,
and this guarantees that the closed-loop system is~$2$-contracting. Recall that in a~$2$-contracting system 
every nonempty omega limit set is a single equilibrium,
 thus ruling out chaotic attractors and even non-trivial limit cycles~\cite{li1995}.  
Fig.~\ref{fig:chaos_closed} depicts the behaviour of the closed-loop system with~$b=0.1$   and~$c=2b-1.1$. The closed-loop system is thus $2$-contracting,
and as expected  
every solution converges to an equilibrium. 
\end{Example}

\section{$k$-positive   systems}
Ref.~\cite{Eyal_k_posi} introduced the notions of~$k$-positive and~$k$-cooperative systems. The LTV~\eqref{Eq:ltv} is called~$k$-positive if~$A^{[k]}(t)$ is
Metzler   for all~$t$. 
For~$k=1$ this reduces to requiring
that~$A(t)$ is Metzler for all~$t$. In this case the system is positive i.e. the flow maps~$\R^n_+$ to~$\R^n_+$ (and also~$\R^n_-:=-\R^n_+$ to~$\R^n_-$)~\cite{farina2000}.
In other words, the flow maps the set of vectors with zero sign variations to itself.

$k$-positive systems map the set of vectors    with up to~$k-1$ sign variations to itself. To explain this, we recall some definitions and results from the theory of totally positive~(TP) matrices, that is, matrices whose minors are all positive~\cite{total_book,pinkus}.

For a vector~$x\in\R^n\setminus\{0\}$, let~$s^-(x) $ denote the number of sign variations in~$x$ after deleting all  its
zero entries. For example,~$s^-(\begin{bmatrix}-1&0&0&2&-3\end{bmatrix}^T)=2$. We define~$s^-(0):=0$. For a vector~$x\in\R^n $,  let~$s^+(x)$
denote the maximal possible number of sign variations in~$x$ after setting  every zero entry in~$x$ to either~$-1$ or~$+1$. For example,~$s^+(\begin{bmatrix}-1&0&0&2&-3\end{bmatrix}^T)=4$. These definitions imply that
$
0\leq s^-(x)\leq s^+(x)\leq n-1$,   for all $ x\in\R^n.
$

For any $k\in \{1,\dots,n\}$,
define the sets 
\begin{align}\label{eq:Pksets}
P_{-}^{k} &:= \{ z\in\mathbb{R}^{n}: \; s^{-} (z) \le k-1 \},\nonumber \\
P_{+}^{k} &:= \{ z\in \mathbb{R}^{n}: \; s^{+} (z) \le k-1 \}.
\end{align}
In other words,
these are the sets of all vectors with up to~$k-1$ sign variations. 
Then~$P_{-}^{k}$ is closed, and it can be shown that~$P_{+}^{k}=\Int(P_{-}^{k})$. 
For example,
\begin{align*}
    P_{-}^{1} = \mathbb{R}_{+}^{n} \cup \mathbb{R}_{-}^{n}, \;\;\; 
    P_{+}^{1} = \Int(\mathbb{R}_{+}^{n}) \cup \Int( \mathbb{R}_{-}^{n}). 
\end{align*}

\begin{Definition}
The LTV~\eqref{Eq:ltv} is called \emph{$k$-positive} on an interval~$[a,b]$ if for any~$a<t_0  < b$,
\[
x(t_0)\in P^k_- \implies x(t,x(t_0)) \in P^k_- \text{ for all }   t_0 \leq t <b ,
\]
and is called \emph{strongly $k$-positive} if
\[
x(t_0)\in P^k_- \implies x(t,x(t_0)) \in P^k_+ \text{ for all }  t_0 < t < b. 
\]
\end{Definition}

In other words, the sets of up to~$k-1$ sign variations are invariant sets of the dynamics.  

An important property   of TP matrices is their sign variation diminishing property: if~$A\in\R^{n\times n}$ is~TP and~$x\in\R^n\setminus\{0\}$ then~$s^+(Ax)\leq s^-(x)$. In other words, multiplying a vector by a TP matrix can only decrease the number of sign variations. For our purposes, we need a more specialized result. Recall that~$A\in\R^{n\times n}$ is called \emph{sign-regular of order~$k$} if   its minors of order~$k$ are all non-positive or all   non-negative, and \emph{strictly 
sign-regular of order~$k$}
if   its minors of order~$k$ are    all positive or all   negative
\begin{Proposition}\label{thm:BenAvraham_SSRk}
%
\cite{CTPDS}
Let $A\in\mathbb{R}^{n\times n}$ be a non-singular matrix. Pick $k\in \{1,\dots,n\}$. Then the following two conditions are equivalent:
\begin{enumerate}
    \item For any $x\in\mathbb{R}^{n}  $ with
    $s^{-}(x) \le k-1$, we have~$s^{-} (Ax) \le k-1$.
    \item $A$ is  sign-regular of order~$k$.
\end{enumerate}
Also,  the following two conditions are equivalent:
\begin{enumerate}[I.]
    \item For any $x\in\mathbb{R}^{n} \setminus \{ 0 \}$ with
    $s^{-}(x) \le k-1$, we have~$s^{+} (Ax) \le k-1$.
    \item $A$ is strictly sign-regular of order~$k$.
\end{enumerate}
\end{Proposition}

Using these tools allows to characterize  
the behaviour of~$k$-positive LTVs.
\begin{Theorem}\label{thm:k_pois_ltv}
The LTV~\eqref{Eq:ltv}
is~$k$-positive on~$[a,b]$ iff~$A^{[k]}(t)$ is Metzler  
for all~$t\in (a,b)$.
It is strongly $k$-positive  on~$[a,b]$
iff~$A^{[k]}(t)$ is Metzler for all~$t\in(a,b)$, 
and~$A^{[k]}(t)$ is irreducible  for all~$t \in(a,b)$  except, perhaps, 
at isolated time points.
\end{Theorem}

The proof is simple. Consider for example  the second assertion in the theorem. The  Metzler and irreducibility assumptions 
imply that the matrix differential system~\eqref{eq:expat} is a positive linear system, and furthermore, that all the entries of~$\Phi^{(k)}(t,t_0)$ are positive for all~$t>t_0$. Thus,~$\Phi(t,t_0)$ is strictly sign-regular of order~$k$ for all~$t>t_0$. Since~$x(t,x(t_0))=\Phi(t,t_0)x(t_0)$,  applying Prop.~\ref{thm:BenAvraham_SSRk} completes the proof. 

This line of reasoning  demonstrates a general and useful principle, namely, given conditions on~$A^{[k]}$ we can   apply standard tools from dynamical  systems theory  to the ``$k$-compound dynamics''~\eqref{eq:expat}, and deduce results on the behaviour of the solution~$x(t)$ of~\eqref{Eq:ltv}. 

A natural question is: when is~$A^{[k]}$ a Metzler matrix? This can be answered using
Prop.~\ref{prop:Explicit_A_k} in terms of   sign pattern conditions on the entries~$a_{ij}$ of~$A$. This is useful as in   fields like
chemistry and systems biology,   exact values of various parameters are typically unknown, but their signs may be inferred from various properties of the system~\cite{sontag_near_2007}.
\begin{Proposition}\label{prop:sign_pattern_for_k_posi}
Let~$A \in\R^{n\times n}$ with~$n\geq 3$.
Then
\begin{enumerate}
\item \label{case:nminus1}
$A^{[ n-1]}$ is Metzler iff
  $a_{ij}\geq 0$ for all~$i,j$ with~$i-j$ odd, and~$a_{ij}\leq 0$ for all~$i,j$ with~$i\not =j$ and~$i-j$ even;
\item \label{case:kodd}
for any  odd~$k  $ in the range~$1<k<n-1$,  
$A^{[k]}$ is Metzler   
 iff
 $a_{1n},a_{n1}\geq 0$, $a_{ij}\geq 0$ for all~$|i-j|=1$, and~$a_{ij}=0$ for all~$1<|i-j|<n-1$;
\item \label{case:keven} 
for
any  even~$k  $ in the range~$1<k<n-1$,  $A^{[k]}$ is Metzler   
 iff
 $a_{1n},a_{n1}\leq 0$, $a_{ij}\geq 0$ for all~$|i-j|=1$, and~$a_{ij}=0$ for all~$1<|i-j|<n-1$.
\end{enumerate}
\end{Proposition}

In Case~\ref{case:nminus1}) there exists a non-singular matrix~$T$ such~$-TAT^{-1}$ is Metzler. In other words, there exists a coordinate transformation such that 
in the new coordinates the dynamics is competitive. 
Thus, $k$-positive systems, with~$k\in\{1,\dots,n-1\}$, may be viewed as a kind of interpolation from cooperative to competitive systems. 
In Case~\ref{case:kodd}),~$A$ is in particular Metzler. Case~\ref{case:keven})  is illustrated in the next example.

\begin{figure}
 \begin{center}
  \includegraphics[scale=0.5]{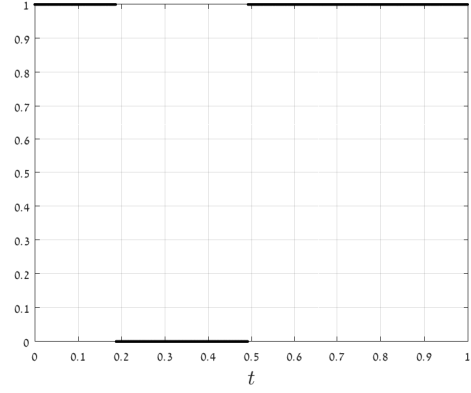}
	\caption{ $s^-(x(t,x(0)))$ as a function of~$t$.}
	\label{fig:signs}
\end{center}
\end{figure}

 \begin{Example}
 Consider the case~$n=3$ and~$A=\begin{bmatrix}
 -1& 1& -2\\ 
 0& 1& 0.1 \\ -3& 0 &1  
 \end{bmatrix} $. Note that~$A$ is not Metzler, yet
$
 A^{[2]}=\begin{bmatrix}
 0& 0.1& 2\\ 0& 0& 1 \\ 3& 0 &2  
 \end{bmatrix}
$
 is Metzler (and irreducible). 
 Thm.~\ref{thm:k_pois_ltv} guarantees that for any~$x(0)$ with~$s^-(x(0))\leq 1$, we have
 \be\label{eq:boundsmin}
 s^-(x(t,x(0)))\leq 1\text{  for all }t\geq 0.
 \ee
 Fig.~\ref{fig:signs}
 depicts~$s^-(x(t,x(0)))=s^-(\exp(At)x(0))$ for~$x(0)=\begin{bmatrix}
 4&-21&-1
 \end{bmatrix}^T$. Note that~$s^-(x(0))=1$. It may be seen  that~$s^-(x(t,x(0)))$ decreases and then increases, but always satisfies the bound~\eqref{eq:boundsmin}.
 \end{Example}

%
%
%
%
         
\subsection{Totally positive differential systems}
A matrix~$A\in\R^{n\times n} $ is called a \emph{Jacobi matrix} if~$A$ is tri-diagonal with positive entries on the super- and sub-diagonals.
An immediate implication of Prop.~\ref{prop:sign_pattern_for_k_posi} is that~$A^{[k]}$ is Metzler and irreducible for all~$k\in\{1,\dots,n-1\}$ iff~$A$ is Jacobi. 
It then follows that for any~$t>0$ 
the matrices~$(\exp (At) )^{(k)} $, $k=1,\dots,n$, are   positive, that is,~$\exp(At)$ is TP for all~$t>0$. 
Combining this with Thm.~\ref{thm:k_pois_ltv} yields the following. 
\begin{Proposition}\cite{schwarz1970}\label{thm:TP}
The following two conditions are equivalent. 
\begin{enumerate} 
\item $A$ is Jacobi.
\item  for any~$x_0\in\R^n \setminus\{0\}$ 
the solution of the LTI~$\dot x(t)=Ax(t)$, $x(0)=x_0$,
satisfies
\[
s^+(x(t, x_0) )\leq s^-(x_0) \text { for all } t>0.
\]
\end{enumerate}
\end{Proposition}

In other words,~$s^-(x(t,x_0))$ 
and also~$s^+(x(t,x_0))$ are non-increasing  functions of~$t$, and may thus be considered as  piece-wise constant Lyapunov functions for the dynamics. 

Prop.~\ref{thm:TP} was proved by Schwarz~\cite{schwarz1970}, yet he  only considered  linear systems. It was recently shown~\cite{margaliot2019revisiting} that important results on the asymptotic behaviour 
of time-invariant and periodic 
time-varying nonlinear systems with a Jacobian that is a Jacobi matrix for all~$t,x$~\cite{smillie,periodic_tridi_smith} follow from the fact that the associated  variational equation is a totally positive~LTV.

\section{$k$-cooperative    systems}
 
%
We now review the applications of~$k$-positivity
to the time-invariant nonlinear system:
\be\label{eq:time_invariant_non_linear}
\dot x=f(x),
\ee
with~$f\in C^1$. Let~$J(x):=\frac{\partial }{\partial x}f(x)$. We assume that the trajectories of~\eqref{eq:time_invariant_non_linear} evolve on a  convex and compact  state-space~$\Omega\subseteq\R^n$.

Recall that~\eqref{eq:time_invariant_non_linear}
is called \emph{cooperative} if~$J(x)$ is Metzler for all~$x\in \Omega$. In other words, the variational equation associated with~\eqref{eq:time_invariant_non_linear}
is positive. The slightly stronger condition of strong cooperativity
has far reaching implications. By Hirsch's quasi-convergence theorem~\cite{hlsmith},
almost every bounded trajectory  
converges to the set of equilibria. 

It is natural to define~$k$-cooperativity by requiring  that the variational equation associated   
with~\eqref{eq:time_invariant_non_linear} is~$k$-positive.
 

\begin{Definition}\label{def:k-coop}\cite{Eyal_k_posi}
The nonlinear system~\eqref{eq:time_invariant_non_linear} is called  \emph{[strongly] $k$-cooperative} 
if the associated  LTV~\eqref{eq:var_eqn} is [strongly]~$k$-positive for any~$a,b\in\Omega$.
\end{Definition}

Note that for~$k=1$ this reduces to the definition of a cooperative  [strongly coopertive] dynamical system.

One immediate implication of Definition~\ref{def:k-coop} is the existence of certain invariant sets of the dynamics.

 
\begin{Proposition} \label{prop:inhgt}
Suppose that~\eqref{eq:time_invariant_non_linear}  
is~$k$-cooperative.
Pick~$a,b\in\Omega$. Then
\[
a-b\in P^k_- \implies  x(t,a)-x(t,b) \in P^k_-  \text{ for all } t \geq  0. 
\]
If, furthermore,~$0 \in \Omega $ and~$0$
is an equilibrium point of~\eqref{eq:time_invariant_non_linear}, i.e.~$f(0)=0$   then 
\[
a \in P^k_- \implies  x(t,a)  \in P^k_-  \text{ for all } t \geq  0. 
\]
\end{Proposition}

The sign pattern conditions in Prop.~\ref{prop:sign_pattern_for_k_posi} can be used to provide simple to verify sufficient conditions for
[strong] $k$-cooperativity  of~\eqref{eq:time_invariant_non_linear}. Indeed, if~$J(x)$ satisfies a sign pattern condition for all~$x\in \Omega$ then   the integral of~$J$ in the variational equation~\eqref{eq:var_eqn} satisfies the same sign pattern, and thus so does~$A^{ab}$. 
The next example, adapted from~\cite{Eyal_k_posi}, 
illustrates this.

\begin{Example} Ref.~\cite{Elkhader1992}
studied the   nonlinear system
\begin{align}\label{eq:alexsys}
\dot x_1&=f_1(x_1,x_n),\nonumber\\
\dot x_i &= f_i(x_{i-1},x_i,x_{i+1}),
\quad i=2,\dots,n-1,\nonumber\\
\dot x_n&=f_n(x_{n-1},x_n),
\end{align}
with the following assumptions:
the state-space~$\Omega\subseteq\R^n$
 is convex, $f_i\in C^{n-1}$, $i=1,\dots,n$, and 
there exist~$\delta_i\in\{-1,1\}$, $i=1,\dots,n$,
 such that
\begin{align*}
\delta_1\frac{\partial }{\partial x_n}f_1(x)   &>0,\\
\delta_2\frac{\partial }{\partial x_1}f_2(x) , \delta_3\frac{\partial }{\partial x_3}f_2(x)&>0,\\
&\vdots\\
\delta_{n-1} \frac{\partial }{\partial x_{n-2}}f_{n-1}(x),
\delta_n \frac{\partial }{\partial x_{n}}f_{n-1}(x)
 &>0,\\
\delta_n \frac{\partial }{\partial x_{n-1}}f_n(x)&>0,
\end{align*}
for all~$x\in\Omega$. 
This is a generalization of     the
monotone cyclic feedback system analyzed  in~\cite{poin_cyclic}.
As noted in~\cite{Elkhader1992}, we may assume without loss of generality that~$\delta_2=\delta_3=\dots=\delta_n=1$ and~$\delta_1 \in \{-1,1\}$. Then the Jacobian of~\eqref{eq:alexsys} has the form
\[
J(x)=\begin{bmatrix}
*& 0 &0 &0 &\dots & 0 & 0  & \sgn(\delta_1) \\
>0 & * &>0 &0 &\dots & 0& 0 & 0  \\
0& >0 & * &>0  &\dots &0&  0 & 0 \\
&&&\vdots\\
0&  0 & 0 & 0  &\dots & 0  &>0&  * \\
\end{bmatrix} ,
\]
for all~$x\in \Omega$. Here~$*$ denotes ``don't care''. Note that~$J(x)$ is irreducible for all~$x\in \Omega$. 

If~$\delta_1=1$ then~$J(x)$ is Metzler, so the system is strongly~$1$-cooperative. 

If~$\delta_1=-1$ then~$J(x)$ satisfies the sign pattern in Case~\ref{case:keven} in Prop.~\ref{prop:sign_pattern_for_k_posi}, so the system is strongly $2$-cooperative.
(If~$n$ is odd then~$J(x)$ also satisfies the sign pattern in Case~\ref{case:nminus1}, so  there is a coordinate transformation for which the  system is also strongly competitive.) 
\end{Example}

 The main result in~\cite{Eyal_k_posi} is that strongly~$2$-cooperative systems satisfy  
 a  
 strong Poincar\'{e}-Bendixson property.
		\begin{Theorem}  \label{thm:2dim}
	 Suppose that~\eqref{eq:time_invariant_non_linear} is strongly~$2$-cooperative. Pick~$a\in \Omega$.
If the omega limit set~$\omega(a)$ does not include an equilibrium 
then it is a closed orbit.
	\end{Theorem}

The proof of this  result is based on the seminal results of Sanchez~\cite{sanchez2009cones}. Yet, it is  considerably
stronger than  the main result in~\cite{sanchez2009cones}, as it applies to \emph{any} trajectory emanating from~$\Omega$   and not only to so called  \emph{pseudo-ordered} trajectories (see the definition in~\cite{sanchez2009cones}). 

The Poincar\'{e}-Bendixon property is useful because often it can be combined with a local analysis near  the equilibrium points to provide a global picture  of the dynamics. For a recent application of Thm.~\ref{thm:2dim} to a model from 
  systems biology, see~\cite{Margaliot868000}.

\section{Conclusion}
$k$-compound matrices describe the evolution of~$k$-dimensional polygons along an LTV dynamics. This geometric property has important consequences in systems and control theory. 
This holds for  both LTVs and also time-varying nonlinear systems, as their variational equation is an LTV. 

Due to space limitations, we  considered here only a partial list of applications. Another application, for example, is based on generalizing  diagonal stability for the LTI~$\dot x= Ax$ to~\emph{$k$-diagonal stability} by requiring that there exists a diagonal  and positive-definite matrix~$D$ such that
$
D A^{[k]}+(A^{[k]})^T D 
$
is negative-definite~\cite{cheng_diag_stab}.

Another interesting line of research is based on analyzing  systems with inputs and outputs. A SISO system is called  \emph{externally~$k$-positive}
if any input with up to~$k$ sign variations induces an output with up to~$k$ sign variations~\cite{gruss1,gruss2,gruss3,gruss4,gruss5}. For LTIs with a zero initial condition the input-output mapping is described by a convolution with the impulse response and then external~$k$-positivity 
is  related to interesting
results in statistics~\cite{ibragimov1956}
and the theory of infinite-dimensional linear operators~\cite{karlin_tp}.


\begin{thebibliography}{10}
\providecommand{\url}[1]{#1}
\csname url@samestyle\endcsname
\providecommand{\newblock}{\relax}
\providecommand{\bibinfo}[2]{#2}
\providecommand{\BIBentrySTDinterwordspacing}{\spaceskip=0pt\relax}
\providecommand{\BIBentryALTinterwordstretchfactor}{4}
\providecommand{\BIBentryALTinterwordspacing}{\spaceskip=\fontdimen2\font plus
\BIBentryALTinterwordstretchfactor\fontdimen3\font minus
  \fontdimen4\font\relax}
\providecommand{\BIBforeignlanguage}[2]{{%
\expandafter\ifx\csname l@#1\endcsname\relax
\typeout{** WARNING: IEEEtranS.bst: No hyphenation pattern has been}%
\typeout{** loaded for the language `#1'. Using the pattern for}%
\typeout{** the default language instead.}%
\else
\language=\csname l@#1\endcsname
\fi
#2}}
\providecommand{\BIBdecl}{\relax}
\BIBdecl

\bibitem{rola_spect}
R.~Alseidi, M.~Margaliot, and J.~Garloff, ``On the spectral properties of
  nonsingular matrices that are strictly sign-regular for some order with
  applications to totally positive discrete-time systems,'' \emph{J. Math.
  Anal. Appl.}, vol. 474, pp. 524--543, 2019.

\bibitem{DT_K_POSI}
R.~{Alseidi}, M.~{Margaliot}, and J.~{Garloff}, ``Discrete-time $k$-positive
  linear systems,'' \emph{IEEE Trans.\ Automat.\ Control}, vol.~66, no.~1, pp.
  399--405, 2021.

\bibitem{sontag_cotraction_tutorial}
Z.~Aminzare and E.~D. Sontag, ``Contraction methods for nonlinear systems: A
  brief introduction and some open problems,'' in \emph{{Proc.\ 53rd IEEE Conf.
  on Decision and Control}}, Los Angeles, CA, 2014, pp. 3835--3847.

\bibitem{chaos_survey}
V.~Basios, C.~G. Antonopoulos, and A.~Latifi, ``Labyrinth chaos: Revisiting the
  elegant, chaotic, and hyperchaotic walks,'' \emph{Chaos}, vol.~30, no.~11, p.
  113129, 2020.

\bibitem{CTPDS}
T.~Ben-Avraham, G.~Sharon, Y.~Zarai, and M.~Margaliot, ``Dynamical systems with
  a cyclic sign variation diminishing property,'' \emph{IEEE Trans.\ Automat.\
  Control}, vol.~65, pp. 941--954, 2020.

\bibitem{coppel1965stability}
W.~A. Coppel, \emph{Stability and {A}symptotic {B}ehavior of {D}ifferential
  {E}quations}.\hskip 1em plus 0.5em minus 0.4em\relax Boston, MA: D. C. Heath,
  1965.

\bibitem{Douady1980}
A.~Douady and J.~Oesterl\'{e}, ``Dimension de {Hausdorff} des attracteurs,''
  \emph{C. R. Acad. Sc. Paris}, vol. 290, pp. 1135--1138, 1980.

\bibitem{Elkhader1992}
A.~S. Elkhader, ``A result on a feedback system of ordinary differential
  equations,'' \emph{J. Dyn. Diff. Equat.}, vol.~4, no.~3, pp. 399--418, 1992.

\bibitem{total_book}
S.~M. Fallat and C.~R. Johnson, \emph{Totally Nonnegative Matrices}.\hskip 1em
  plus 0.5em minus 0.4em\relax Princeton, NJ: Princeton University Press, 2011.

\bibitem{farina2000}
L.~Farina and S.~Rinaldi, \emph{Positive Linear Systems: Theory and
  Applications}.\hskip 1em plus 0.5em minus 0.4em\relax John Wiley, 2000.

\bibitem{fiedler_book}
M.~Fiedler, \emph{Special Matrices and Their Applications in Numerical
  Mathematics}, 2nd~ed.\hskip 1em plus 0.5em minus 0.4em\relax Mineola, NY:
  Dover Publications, 2008.

\bibitem{notes_comb_algebra}
\BIBentryALTinterwordspacing
D.~Grinberg, ``Notes on the combinatorial fundamentals of algebra,'' 2020.
  [Online]. Available: \url{https://arxiv.org/abs/2008.09862}
\BIBentrySTDinterwordspacing

\bibitem{gruss1}
C.~Grussler, T.~B. Burghi, and S.~Sojoudi, ``Internally {H}ankel-positive
  systems,'' 2021, arXiv preprint arXiv:2103.06962.

\bibitem{gruss3}
C.~Grussler, T.~Damm, and R.~Sepulchre, ``Balanced truncation of $k$-positive
  systems,'' 2020, arXiv preprint arXiv:2006.13333.

\bibitem{gruss5}
C.~Grussler and R.~Sepulchre, ``Strongly unimodal systems,'' in \emph{Proc.
  18th Euro. Control Conf.}, 2019, pp. 3273--3278.

\bibitem{gruss2}
------, ``Variation diminishing {H}ankel operators,'' in \emph{{Proc.\ 59th
  IEEE Conf. on Decision and Control}}, 2020, pp. 4529--4534.

\bibitem{gruss4}
------, ``Variation diminishing linear time-invariant systems,'' 2020, arXiv
  preprint arXiv:2006.10030.

\bibitem{ibragimov1956}
I.~Ibragimov, ``On the composition of unimodal distributions,'' \emph{Theory of
  Probability \& Its Applications}, vol.~1, no.~2, p. 255–260, 1956.

\bibitem{karlin_tp}
S.~Karlin, \emph{Total Positivity, Volume 1}.\hskip 1em plus 0.5em minus
  0.4em\relax Stanford, CA: Stanford University Press, 1968.

\bibitem{diag_Stab_book}
E.~Kaszkurewicz and A.~Bhaya, \emph{Matrix Diagonal Stability in Systems and
  Computation}.\hskip 1em plus 0.5em minus 0.4em\relax New York, NY: Springer,
  2000.

\bibitem{rami_osci}
R.~Katz, M.~Margaliot, and E.~Fridman, ``Entrainment to subharmonic
  trajectories in oscillatory discrete-time systems,'' \emph{Automatica}, vol.
  116, p. 108919, 2020.

\bibitem{SEIR_LI_MULD1995}
M.~Y. Li and J.~S. Muldowney, ``Global stability for the {SEIR} model in
  epidemiology,'' \emph{Math. Biosciences}, vol. 125, no.~2, pp. 155--164,
  1995.

\bibitem{li1995}
------, ``On {R. A. Smith's} autonomous convergence theorem,'' \emph{Rocky
  Mountain J. Math.}, vol.~25, no.~1, pp. 365--378, 1995.

\bibitem{LOHMILLER1998683}
W.~Lohmiller and J.-J.~E. Slotine, ``On contraction analysis for non-linear
  systems,'' \emph{Automatica}, vol.~34, pp. 683--696, 1998.

\bibitem{poin_cyclic}
J.~Mallet-Paret and H.~L. Smith, ``The {P}oincar{\'e}-{B}endixson theorem for
  monotone cyclic feedback systems,'' \emph{J. Dyn. Differ. Equ.}, vol.~2,
  no.~4, pp. 367--421, 1990.

\bibitem{weak_manchester}
\BIBentryALTinterwordspacing
I.~R. Manchester and J.-J.~E. Slotine, ``Combination properties of weakly
  contracting systems,'' 2014. [Online]. Available:
  \url{https://arxiv.org/abs/1408.5174}
\BIBentrySTDinterwordspacing

\bibitem{Margaliot868000}
\BIBentryALTinterwordspacing
M.~Margaliot and E.~D. Sontag, ``Compact attractors of an antithetic integral
  feedback system have a simple structure,'' 2019. [Online]. Available:
  \url{https://www.biorxiv.org/content/early/2019/12/08/868000}
\BIBentrySTDinterwordspacing

\bibitem{margaliot2019revisiting}
------, ``Revisiting totally positive differential systems: A tutorial and new
  results,'' \emph{Automatica}, vol. 101, pp. 1--14, 2019.

\bibitem{muldowney1990}
J.~S. Muldowney, ``Compound matrices and ordinary differential equations,''
  \emph{Rocky Mountain J. Math.}, vol.~20, no.~4, pp. 857--872, 12 1990.

\bibitem{searial12contracting}
R.~Ofir, M.~Margaliot, Y.~Levron, and J.-J. Slotine, ``Serial interconnections
  of~$1$-contracting and~$2$-contracting~systems,'' 2021, submitted.

\bibitem{pinkus}
A.~Pinkus, \emph{Totally Positive Matrices}.\hskip 1em plus 0.5em minus
  0.4em\relax Cambridge, UK: Cambridge University Press, 2010.

\bibitem{sanchez2009cones}
L.~A. Sanchez, ``Cones of rank 2 and the {P}oincar{\'e}-{B}endixson property
  for a new class of monotone systems,'' \emph{J. Diff. Eqns.}, vol. 246,
  no.~5, pp. 1978--1990, 2009.

\bibitem{schwarz1970}
B.~Schwarz, ``Totally positive differential systems,'' \emph{Pacific J. Math.},
  vol.~32, no.~1, pp. 203--229, 1970.

\bibitem{smillie}
J.~Smillie, ``Competitive and cooperative tridiagonal systems of differential
  equations,'' \emph{SIAM J. Math. Anal.}, vol.~15, pp. 530--534, 1984.

\bibitem{periodic_tridi_smith}
H.~L. Smith, ``Periodic tridiagonal competitive and cooperative systems of
  differential equations,'' \emph{SIAM J. Math. Anal.}, vol.~22, no.~4, pp.
  1102--1109, 1991.

\bibitem{hlsmith}
------, \emph{Monotone Dynamical Systems: An Introduction to the Theory of
  Competitive and Cooperative Systems}, ser. Mathematical Surveys and
  Monographs.\hskip 1em plus 0.5em minus 0.4em\relax Providence, RI: Amer.
  Math. Soc., 1995, vol.~41.

\bibitem{sontag_near_2007}
E.~D. Sontag, ``Monotone and near-monotone biochemical networks,''
  \emph{Systems and Synthetic Biology}, vol.~1, pp. 59--87, 2007.

\bibitem{thomas99}
R.~Thomas, ``Deterministic chaos seen in terms of feedback circuits: Analysis,
  synthesis, ``labyrinth chaos'','' \emph{Int. J. Bifurc. Chaos.}, vol.~9,
  no.~10, pp. 1889--1905, 1999.

\bibitem{vid}
M.~Vidyasagar, \emph{Nonlinear Systems Analysis}.\hskip 1em plus 0.5em minus
  0.4em\relax Englewood Cliffs, NJ: Prentice Hall, 1978.

\bibitem{Eyal_k_posi}
E.~Weiss and M.~Margaliot, ``A generalization of linear positive systems with
  applications to nonlinear systems: Invariant sets and the
  {Poincar\'{e}-Bendixson} property,'' \emph{Automatica}, vol. 123, p. 109358,
  2021.

\bibitem{9107214}
E.~{Weiss} and M.~{Margaliot}, ``Is my system of odes k-cooperative?''
  \emph{IEEE Control Systems Letters}, vol.~5, no.~1, pp. 73--78, 2021.

\bibitem{kordercont}
\BIBentryALTinterwordspacing
C.~Wu, I.~Kanevskiy, and M.~Margaliot, ``$k$-order contraction: theory and
  applications,'' 2020, submitted. [Online]. Available:
  \url{https://arxiv.org/abs/2008.10321}
\BIBentrySTDinterwordspacing

\bibitem{cheng_diag_stab}
\BIBentryALTinterwordspacing
C.~Wu and M.~Margaliot, ``Diagonal stability of discrete-time $k$-positive
  linear systems with applications to nonlinear systems,'' 2020, submitted.
  [Online]. Available: \url{arXiv:2102.02144}
\BIBentrySTDinterwordspacing

\bibitem{pines2021}
C.~Wu, R.~Pines, M.~Margaliot, and J.-J. Slotine, ``Generalization of the
  multiplicative and additive compounds of square matrices and contraction in
  the {Hausdorff} dimension,'' 2020, arXiv 2012.13441.

\end{thebibliography}
\end{document}